\documentclass[12pt,twoside]{article}
\usepackage{amssymb}

\font\teneufm=eufm10 scaled \magstep1
\font\seveneufm=eufm7 scaled \magstep1
\font\fiveeufm=eufm5  scaled \magstep1
\newfam\eufmfam

\textfont\eufmfam=\teneufm
\scriptfont\eufmfam=\seveneufm
\scriptscriptfont\eufmfam=\fiveeufm
\def\frak#1{{\fam\eufmfam\relax#1}}

\newfam\msbfam
\font\tenmsb=msbm10 scaled \magstep1  \textfont\msbfam=\tenmsb
\font\sevenmsb=msbm7 scaled \magstep1 \scriptfont\msbfam=\sevenmsb
\font\fivemsb=msbm5 scaled \magstep1  \scriptscriptfont\msbfam=\fivemsb

\makeatletter
\def\blfootnote{\xdef\@thefnmark{}\@footnotetext}
\makeatother

\def\dd#1{\raise1.5pt\hbox{$\,\partial\!$}/\raise-2.5pt\hbox{$\!\partial#1\,$}}

\def\tilde{\widetilde}
\def\hat{\widehat}
\def\5#1{{\mathcal #1}}

\def\CC{{\mathbb C}}

\def\FF{{\mathbb F}}

\def\ra{\rightarrow}

\def\GL{\mathop{\rm GL}\nolimits}
\def\SL{\mathop{\rm SL}\nolimits}

\def\Ann{\mathop{\rm Ann}\nolimits}

\def\emb{\mathop{\rm emb}\nolimits}

 \def\HollowBoxx #1#2#3{{\dimen0=#1 \advance\dimen0 by -#2
       \dimen1=#1 \advance\dimen1 by #3
        \vrule height 0pt depth #3 width #2
       \hskip -#3
       \vrule height #1 depth #3 width #3}}
 \def\LeftContraction{\mathord{\kern1.45pt \HollowBoxx{6pt}{3.5pt}{.4pt}}\,}

 \def\HollowBox #1#2#3{{\dimen0=#1 \advance\dimen0 by -#3
       \dimen1=#1 \advance\dimen1 by #3
        \vrule height #1 depth #3 width #3
        \vrule height 0pt depth #3 width #2
        \hskip -#3}}
 \def\RightContraction{\mathord{\, \HollowBox{6pt}{3.1pt}{.4pt}} \kern1.6pt}

\def\qed{{\hfill $\Box$}}
\newtheorem{theorem}{THEOREM}[section]

\newtheorem{example}[theorem]{Example}
\newtheorem{remark}[theorem]{Remark}
\newtheorem{proposition}[theorem]{Proposition}
\newtheorem{conjecture}[theorem]{Conjecture}

\begin{document}

\begin{center}
{\Large \bf Invariants of Artinian Gorenstein Algebras
\vspace{0.3cm}\\
and Isolated Hypersurface Singularities}\\
\vspace{0.5cm}
\normalsize M. G. Eastwood and A. V. Isaev
\end{center}

\begin{quotation} 
{\small \sl \noindent We survey our recently proposed method for constructing biholomorphic invariants of quasihomogeneous isolated hypersurface singularities and, more generally, invariants of graded Artinian Gorenstein algebras. The method utilizes certain polynomials associated to such algebras, called nil-polynomials, and we compare them with two other classes of polynomials that have also been used to produce invariants.}
\end{quotation}

\thispagestyle{empty}

\pagestyle{myheadings}
\markboth{M. G. Eastwood and A. V. Isaev}{Invariants of Artinian Gorenstein Algebras}

\setcounter{section}{0}

\section{Introduction}\label{intro}
\setcounter{equation}{0}

On $27^{\mathrm{th}}$ October 2001, one of us (MGE) gave a talk entitled
`Invariants of isolated hypersurface singularities' at the West Coast Lie
Theory Workshop held at the University of California, Berkeley. The purpose of
this talk was to propose, at that time only by means of examples, a method for
extracting invariants of isolated hypersurface singularities. Since classical
invariant theory was the key ingredient in the method and since classical
invariants may be derived from representation theory of the general linear
group, this seemed an appropriate topic for the West Coast Lie Theory Series.
The examples from this talk were presented in~\cite{Ea}. Since that time, the
method from~\cite{Ea} has been considerably improved (by authors other than
MGE) and related to other methods. Also, it has been realized that the proper
framework in which to formulate the results is within the theory of Gorenstein
algebras. Although this now seems quite far from Lie theory, we believe it is
the natural development and we thank the editors of this volume for the
opportunity to present this survey, which might otherwise appear out of place.

Let $V$ be a complex hypersurface germ with an isolated singularity at $0\in\CC^n$, $n\ge 2$. It follows that the complex structure of $V$ is reduced, i.e.~defined by the ideal $I(V)$ in the algebra ${\mathcal O}_n$ of germs of holomorphic functions at the origin that consists of all elements of  ${\mathcal O}_n$ vanishing on $V$. The singularity of $V$ is called {\it quasihomogeneous}\, if some (hence every) generator of $I(V)$ in some coordinates $z_1,\dots,z_n$ near the origin is the germ of a quasihomogeneous polynomial $Q(z_1,\dots,z_n)$, i.e.~a polynomial satisfying $Q(t^{p_1}z_1,\dots,t^{p_n}z_n)\equiv t^qQ(z_1,\dots,z_n)$ for fixed positive integers $p_1,\dots,p_n,q$ and all $t\in\CC$, where $p_j$ is called the weight of $z_j$ and $q$ the weight of $Q$. The singularity is called {\it homogeneous}\, if one can choose $Q$ to be a form (homogeneous polynomial). This paper concerns the biholomorphic invariants of quasihomogeneous singularities introduced in our recent article \cite{EI}. Various invariants of hypersurface singularities have been extensively studied by many authors (see Chapter 1 in \cite{GLS} for an account of some of the results). The objective of \cite{EI} was to construct, in the quasihomogeneous case, numerical invariants that: (i) are easy to compute, and (ii) form a complete system with respect to the biholomorphic equivalence problem for hypersurface germs. 
Although we succed only in a very limited setting (concerning binary quintics 
and sextics, whose classical invariants are well-understood),
to the best of our knowledge, no such system of invariants had been previously
known.

Our approach utilizes the {\it Milnor algebra}\, of $V$, which is the complex local commutative associative algebra $M(V):={\mathcal O}_n/J(f)$, where $f$ is any generator of $I(V)$ and $J(f)$ the ideal in ${\mathcal O}_n$ generated by all first-order partial derivatives of $f$ calculated with respect to some coordinate system near the origin. It is easy to observe that the above definition is independent of the choice of $f$ as well as the coordinate system, and that the Milnor algebras of biholomorphically equivalent singularities are isomorphic. Furthermore, the dimension $n$ and the isomorphism class of $M(V)$ determine $V$ up to biholomorphic equivalence (see \cite{Sh} and a more general result in \cite{MY}). Thus, any quantity that depends only on the isomorphism class of $M(V)$ is a biholomorphic invariant of $V$, and any collection of quantities of this kind uniquely 
characterizing 
the isomorphism class of every Milnor algebra is a complete system of biholomorphic invariants for hypersurface germs of fixed dimension. 

In order to produce invariants of Milnor algebras of quasihomogeneous singularities, we focus on three important properties of $M(V)$. First, since the singularity of $V$ is isolated, one has $\dim_{\CC}M(V)<\infty$ (see, e.g.~Chapter 1 in \cite{GLS}), that is, the algebra $M(V)$ is {\it Artinian}. It then follows that the first-order partial derivatives of $f$ form a regular sequence in ${\mathcal O}_n$ (see Theorem 2.1.2 in \cite{BH}), hence, by \cite{B}, the algebra $M(V)$ is {\it Gorenstein}. (Recall that a local commutative associative algebra $A$ of finite vector space dimension greater than 1 is Gorenstein if the annihilator $\Ann({\mathfrak m}):=\{x\in{\mathfrak m}:x{\mathfrak m}=0\}$ of the maximal ideal ${\mathfrak m}\subset A$ is 1-dimensional -- see, e.g.~\cite{Hu}.) Finally, $M(V)$ is {\it (non-negatively) graded}, i.e.~it can be represented as a direct sum $M(V)=\bigoplus_{i\ge 0}L_i$, where $L_i$ are subspaces such that $L_0\simeq\CC$ and $L_iL_j\subset L_{i+j}$ for all $i,j$. Indeed, choosing coordinates near the origin in which $f$ is the germ of a quasihomogeneous polynomial, we set $L_i$ to be the subspace of elements of $M(V)$ represented by germs of quasi-homogeneous polynomials of weight $i$.

Rather than focussing on invariants of Milnor algebras of quasihomogeneous singularities, one can take a broader viewpoint and introduce certain invariants of general complex graded Artinian Gorenstein algebras. Our method for constructing invariants is based on associating to every algebra $A$ a form $P_A$ of degree $d_A$ on a complex vector space $W_A$ of dimension $N_A$, such that for any pair of isomorphic algebras $A$, $\tilde A$ there exists a linear isomorphism $\varphi:W_A\ra W_{\tilde A}$ with $P_A=P_{\tilde A}\circ\varphi$. Then, upon identification of $W_A$ with  $\CC^{N_A}$, for any absolute classical invariant ${\bf I}$ of forms of degree $d_A$ on $\CC^{N_A}$, the quantity ${\bf I}(P_A)$ is invariantly defined. Observe that for a given choice of $P_A$ all invariants of this kind are easy to calculate using computer algebra.  

The idea of building invariants by the above method goes back at least as far as article \cite{Em} (see also \cite{ER} and references therein for more detail), where it was briefly noted for the case of {\it standard graded}\, Artinian Gorenstein algebras, in which case $P_A$ is one of {\it (Macaulay's) inverse systems}\, for $A$. For instance, if the singularity of $V$ is homogeneous, the algebra $M(V)$ is standard graded. For general quasihomogeneous singularities the idea was explored in \cite{Ea}, with $P_A$ being a certain form $a$ on ${\mathfrak m}/{\mathfrak m}^2$. Finally, in \cite{EI} we utilized homogeneous components of {\it nil-polynomials}\, introduced in \cite{FIKK} to construct a large number of invariants of arbitrary graded Artinian Gorenstein algebras. Relationships among the above three choices of $P_A$ are discussed in Section \ref{sect1}, where we will see, in particular, that nil-polynomials can be regarded as certain extensions of both inverse systems and the form $a$. Hence, the invariants produced from nil-polynomials incorporate those arising from the other two possibilities for $P_A$. The construction of this most general system of invariants and results concerning its completeness are surveyed in Section \ref{sect2}.     

{\bf Acknowledgement.} This work is supported by the Australian Research Council.

\section{Polynomials associated to Artinian\\ Gorenstein algebras}\label{sect1}\label{comparison}
\setcounter{equation}{0}

In this section we establish relationships among three kinds of polynomials arising from Artinian Gorenstein algebras. As mentioned in the introduction, we will consider inverse systems, the form $a$ introduced in \cite{Ea}, and nil-polynomials introduced in \cite{FIKK}. For expository purposes, it is convenient for us to start with nil-polynomials.

Let $A$ be an Artinian Gorenstein algebra over a field $\FF$ of characteristic zero, with $\dim_{\FF}A>2$ and maximal ideal ${\mathfrak m}$. Define a map $\exp: {\mathfrak m}\ra {\bf 1}+{\mathfrak m}$ by the formula
$$
\displaystyle\exp(x):=\sum_{s=0}^{\infty}\frac{1}{s!}x^s,
$$
where $x^0:={\bf 1}$, with ${\bf 1}$ being the identity element of $A$. By Nakayama's lemma, ${\mathfrak m}$ is a nilpotent algebra, and therefore the above sum is in fact finite, with the highest-order term corresponding to $s=\nu$, where $\nu\ge 2$ is the nil-index of ${\mathfrak m}$ (i.e.~the largest of all integers $\mu$ for which ${\frak m}^{\mu}\ne 0$). Fix a hyperplane $\Pi$ in ${\mathfrak m}$ complementary to $\Ann({\mathfrak m})={\frak m}^{\nu}$. An $\FF$-valued polynomial $P$ on $\Pi$ is called a nil-polynomial if there exists a linear form $\omega:A\ra\FF$ such that $\ker\omega=\langle\Pi,{\bf 1}\rangle$ and
$$
P=\omega\circ\exp|_{\Pi},\,\,\hbox{i.e.,}\,\, P(x)=\omega\left(\sum_{s=2}^{\nu}\frac{1}{s!}x^s\right),\,\, x\in\Pi,
$$
where $\langle\,\cdot\,\rangle$ denotes linear span.

As shown in \cite{FIKK}, \cite{FK}, \cite{I}, nil-polynomials solve the isomorphism problem for Artinian Gorenstein algebras as follows: $A$ and $\tilde A$ are isomorphic if and only if the graphs $\Gamma\subset \Pi\times\FF$, $\tilde\Gamma\subset \tilde\Pi\times\FF$ of any nil-polynomials $P$, $\tilde P$ arising from $A$, $\tilde A$, respectively, are affinely equivalent, that is, there exists a bijective affine map $\psi:\Pi\times\FF\ra\tilde\Pi\times\FF$ such that $\psi(\Gamma)=\tilde\Gamma$. Furthermore, if both $A$ and $\tilde A$ are graded, then $\Gamma$ and $\tilde\Gamma$ are affinely equivalent if and only if $P$ and $\tilde P$ are linearly equivalent up to scale, i.e.~there exist $c\in\FF^*$ and a linear isomorphism $\varphi:\Pi\ra\tilde\Pi$ with $cP=\tilde P\circ\varphi$. As will be seen in Section \ref{sect2}, this is exactly the property that allows one to use nil-polynomials for producing invariants of graded Artinian Gorenstein algebras. 

Further, any nil-polynomial $P=\omega\hspace{0.01cm}\circ\hspace{0.01cm}\exp|_{\Pi}$ arising from a Gorenstein algebra $A$ extends to the polynomial $\hat P:=\omega\circ\exp$ on all of ${\mathfrak m}$. Let 
$$
\hat P^{[s]}(x):=\frac{1}{s!}\omega(x^s),\,\, x\in{\mathfrak m},
$$
 be the homogeneous component of $\hat P$ of degree $s$, with $s=2,\dots,\nu$. One has $\hat P^{[s]}(y)=0$, $\hat P^{[s]}(x+y)=\hat P^{[s]}(x)$ for all $x\in{\mathfrak m}$, $y\in{\mathfrak m}^{\nu+2-s}$. Thus, $\hat P^{[s]}$ gives rise to a form ${\bf P}^{[s]}$ on the quotient ${\mathfrak m}/{\mathfrak m}^{\nu+2-s}$. The forms ${\bf P}^{[s]}$ will be used in Section \ref{sect2} for constructing the invariants mentioned above. Here we only observe that the highest-degree form ${\bf P}^{[\nu]}$ defined on ${\mathfrak m}/{\mathfrak m}^2$ is special. Indeed, for any other nil-polynomial $P'$ arising from $A$ the corresponding form ${\bf P}^{'[\nu]}$ coincides with ${\bf P}^{[\nu]}$ up to scale. Moreover, it is clear from the definition of the form $a$ given on p.~305 in \cite{Ea} that it is equal, up to scale, to ${\bf P}^{[\nu]}$. Thus, loosely speaking, $a$ can be regarded, up to proportionality, as the highest-degree homogeneous component of any nil-polynomial.

Next, let $k:=\emb\dim A:=\dim_{\FF}{\mathfrak m}/{\mathfrak m}^2\ge 1$ be the embedding dimension of $A$. Since $\dim_{\FF}A>2$, the hyperplane $\Pi$ contains a $k$-dimensional subspace that forms a complement to ${\mathfrak m}^2$ in ${\mathfrak m}$. Fix any such subspace $L$, choose a basis $e_1,\dots,e_k$ in it, and let $y_1,\dots,y_k$ be the coordinates with respect to this basis. Denote by $R\in\FF[y_1,\dots,y_k]$ the restriction of the nil-polynomial $P$ to $L$ expressed in these coordinates. Clearly, one has
$$
R(y_1,\dots,y_k)=\sum_{j=0}^{\nu}\frac{1}{j!}\omega\Bigl((y_1e_1+\dots+y_ke_k)^j\Bigr),
$$
and the homogeneous component $R^{[\nu]}$ of degree $\nu$ of $R$ is given by
$$
R^{[\nu]}(y_1,\dots,y_k)=\frac{1}{\nu!}\omega\Bigl((y_1e_1+\dots+y_ke_k)^{\nu}\Bigr).
$$
Thus, identifying $L$ with ${\mathfrak m}/{\mathfrak m}^2$, we see that $R^{[\nu]}$ is a coordinate representation of the form ${\bf P}^{[\nu]}$ and therefore that of the form $a$ up to a scaling factor.

Further, the elements $e_1,\dots,e_k$ generate $A$ as an algebra, hence $A$ is isomorphic to $\FF[x_1,\dots,x_k]/I$, where $I$ is the ideal of all relations among $e_1,\dots,e_k$, i.e.~polynomials $f\in\FF[x_1,\dots,x_k]$ with $f(e_1,\dots,e_k)=0$. Observe that $I$ contains the monomials $x_1^{\nu+1},\dots,x_k^{\nu+1}$, and therefore $A$ is also isomorphic to $\FF[[x_1,\dots,x_k]]/\FF[[x_1,\dots,x_k]]I$. It is well-known that, since the quotient $\FF[x_1,\dots,x_k]/I$ is Gorenstein, there is a polynomial $g\in\FF[y_1,\dots,y_k]$ of degree $\nu$ satisfying $\Ann(g)=I$, where
$$
\Ann(g):=\left\{f\in\FF[x_1,\dots,x_k]: f\left(\frac{\partial}{\partial y_1},\dots,\frac{\partial}{\partial y_k}\right)(g)=0\right\}
$$
is the annihilator of $g$ (see, e.g.~\cite{ER} and references therein). The freedom in choosing $g$ with $\Ann(g)=I$ is fully understood, and any such polynomial is called (Macaulay's) inverse system for the Gorenstein quotient $\FF[x_1,\dots,x_k]/I$. The classical correspondence $I\leftrightarrow g$ can be also derived from the {\it Matlis duality}\, (see Section 5.4 in \cite{SV}).
   
Inverse systems can be used for solving the isomorphism problem for quotients of this kind. Namely, two Gorenstein quotients are isomorphic if and only if their inverse systems are equivalent in a certain sense (see Proposition 16 in \cite{Em} and a more explicit formulation in Proposition 2.2 in \cite{ER}). Observe that in general the equivalence relation for inverse systems is harder to analyze than the affine equivalence of graphs of nil-polynomials mentioned above, and therefore the criterion for isomorphism of Artinian Gorenstein algebras in terms of inverse systems seems to be less convenient in applications than that in terms of nil-polynomials. There is one case, however, when the criterion in terms of inverse systems is rather useful. It is discussed in Remark \ref{stangraded} at the end of this section. 

The following theorem provides a connection between nil-polynomials and inverse systems.

\begin{theorem}\label{main}\sl The polynomial $R$ is an inverse system for the quotient $\FF[x_1,\dots,x_k]/I$.
\end{theorem}

\noindent{\bf Proof:} Fix any polynomial $f\in\FF[x_1,\dots,x_k]$
$$
f=\sum_{0\le i_1,\dots,i_k\le N}a_{i_1,\dots,i_k}x_1^{i_1}\dots x_k^{i_k}
$$ 
and calculate
$$
\displaystyle f\left(\frac{\partial}{\partial y_1},\dots,\frac{\partial}{\partial y_k}\right)(R)=\sum_{0\le i_1,\dots,i_k\le N}a_{i_1,\dots,i_k}\sum_{j= i_1+\dots+i_k}^{\nu}\frac{1}{(j-(i_1+\dots+i_k))!}\times
$$
\begin{equation}
\begin{array}{l}
\displaystyle\omega\Bigl((y_1e_1+\dots+y_ke_k)^{j-(i_1+\dots+i_k)}e_1^{i_1}\dots e_k^{i_k}\Bigr)=\\
\vspace{-0.3cm}\\
\displaystyle\sum_{m=0}^{\nu}\frac{1}{m!}\omega\Bigl((y_1e_1+\dots+y_ke_k)^m
\sum_{\mbox{\tiny$\begin{array}{l} 0\le i_1,\dots,i_k\le N,\\\vspace{-0.1cm}\\
i_1+\dots+i_k\le\nu-m\end{array}$}}
a_{i_1,\dots,i_k}e_1^{i_1}\dots e_k^{i_k}\Bigr)=\\
\vspace{-0.3cm}\\
\displaystyle\sum_{m=0}^{\nu}\frac{1}{m!}\omega\Bigl((y_1e_1+\dots+y_ke_k)^m\,f(e_1,\dots,e_k)\Bigr).
\end{array}\label{diff}
\end{equation}
Formula (\ref{diff}) immediately implies $I\subset\Ann(R)$.

Conversely, let $f\in\FF[x_1,\dots,x_k]$ be an element of $\Ann(R)$. Then (\ref{diff}) yields
\begin{equation}
\sum_{m=0}^{\nu}\frac{1}{m!}\omega\Bigl((y_1e_1+\dots+y_ke_k)^m\,f(e_1,\dots,e_k)\Bigr)=0.\label{eq2}
\end{equation}
Collecting the terms containing $y_1^{i_1}\dots y_k^{i_k}$ in (\ref{eq2}) we obtain
\begin{equation}
\omega\Bigl(e_1^{i_1}\dots e_k^{i_k}\,f(e_1,\dots,e_k)\Bigr)=0\label{eq3}
\end{equation}
for all indices $i_1,\dots,i_k$. Since $e_1,\dots,e_k$ generate $A$, identities (\ref{eq3}) yield 
\begin{equation}
\omega\Bigl(A\, f(e_1,\dots,e_k)\Bigr)=0.\label{nondeg}
\end{equation}
Further, since the bilinear form $(a,b)\mapsto \omega(ab)$ is non-degenerate on $A$ (see, e.g.~p.~11 in \cite{He}), identity (\ref{nondeg}) implies $f(e_1,\dots,e_k)=0$. Therefore $f\in I$, which shows that $I=\Ann(R)$ as required.\qed
\vspace{0.3cm}

\begin{remark}\label{remembdim}\rm Theorem \ref{main} easily generalizes to the case of Artinian Gorenstein quotients $\FF[x_1,\dots,x_m]/I$, where $m$ is not necessarily the embedding dimension of the quotient. Indeed, let $e_1,\dots,e_m$ be the elements of $\FF[x_1,\dots,x_m]/I$ represented by $x_1,\dots,x_m$, respectively, and consider the polynomial in $\FF[y_1,\dots,y_m]$ defined as follows:
$$
S(y_1,\dots,y_m):=\sum_{j=0}^{\nu}\frac{1}{j!}\omega\Bigl((y_1e_1+\dots+y_me_m)^j\Bigr),
$$
where $\omega$ is a linear form on $\FF[x_1,\dots,x_m]/I$ with kernel complementary to $\Ann({\mathfrak m})$. Then arguing as in the proof of Theorem \ref{main}, we see that $S$ is an inverse system for $\FF[x_1,\dots,x_m]/I$. However, if $m>\emb\dim\FF[x_1,\dots,x_m]/I$, this inverse system does not come from restricting a nil-polynomial to a subspace of ${\mathfrak m}$ complementary to ${\mathfrak m}^2$.
\end{remark}

We will now give an example illustrating the relationships among nil-polynomials, inverse systems and the form $a$ established above.

\begin{example}\label{ex}\rm Consider the following one-parameter family of algebras:
$$
A_t:=\FF[x_1,x_2]/(2x_1^3+tx_1x_2^3,tx_1^2x_2^2+2x_2^5),\quad t\in\FF,\,\, t\ne\pm 2.
$$
It is straightforward to verify that every $A_t$ is a Gorenstein algebra of dimension 15 with $\nu=7$ for $t\ne 0$ and $\nu=6$ for $t=0$. 

Consider the following monomials in $\FF[x_1,x_2]$:
$$
x_1,\,\,x_2,\,\,x_1^2,\,\,x_1x_2,\,\,x_2^2,\,\,x_1^2x_2,\,\,x_1x_2^2,\,\,x_2^3,\,\,x_1x_2^3,\,\,x_1^2x_2^2,\,\,x_2^4,\,\,x_1^2x_2^3,\,\,x_1x_2^4,\,\,x_1^2x_2^4.
$$
Let $e_1,\dots,e_{14}$, respectively, be the elements of $A_t$ represented by these monomials. They form a basis of ${\mathfrak m}$. The hyperplane $\Pi:=\langle e_1,\dots,e_{13}\rangle$ in ${\mathfrak m}$ is complementary to $\Ann({\mathfrak m})=\langle e_{14}\rangle$, and we denote by $y_1,\dots,y_{13}$ the coordinates in $\Pi$ with respect to $e_1,\dots,e_{13}$. Further, define $\omega:A_t\ra\FF$ to be the linear form such that $\ker\omega=\langle\Pi,{\bf 1}\rangle$ and $\omega(e_{14})=1$. Then the nil-polynomial $P:=\omega\circ\exp|_{\Pi}$ is expressed in the coordinates $y_1,\dots,y_{13}$ as follows:   
$$
\begin{array}{l}
\displaystyle P(y_1,\dots,y_{13})=\frac{t}{10080}y_2^7-\frac{1}{48}y_2^4\left(y_1^2-\frac{t}{5}y_2y_5\right)+\frac{t}{48}y_1^4y_2-\frac{1}{4}y_1^2y_2^2y_5-\\
\vspace{-0.1cm}\\
\displaystyle\hspace{3cm}\frac{1}{6}y_1y_2^3y_4+\frac{t}{24}y_2^3y_5^2+\frac{t}{48}y_2^4y_8-\frac{1}{24}y_2^4y_3+\hbox{terms of $\deg\le 4\,$}.
\end{array}
$$

Further, setting $L:=\langle e_1,e_2\rangle$ and restricting $P$ to $L$, we arrive at an inverse system of $A_t$:
$$
R(y_1,y_2)=\frac{t}{10080}y_2^7-\frac{1}{48}y_1^2y_2^4+\frac{t}{48}y_1^4y_2
$$
(note that all terms of $\deg\le 4\,$ in $P$ vanish for $y_3=\dots=y_{13}=0$). Then, choosing the highest-degree terms in $R$ and identifying $L$ with ${\mathfrak m}/{\mathfrak m}^2$, we obtain the following coordinate representation of the form $a$:
$$
a(y_1,y_2)=\left\{
\begin{array}{ll}
\displaystyle\frac{t}{10080}y_2^7 & \hbox{if $t\ne 0$},\\
\vspace{-0.1cm}\\
\displaystyle-\frac{1}{48}y_1^2y_2^4 & \hbox{if $t=0$}
\end{array}
\right.\quad\hbox{up to scale.}
$$

\end{example}
\vspace{0.3cm}

\begin{remark}\label{stangraded} \rm Let now $A$ be a standard graded algebra, i.e.~an algebra that can be represented as a direct sum 
$A=\bigoplus_{i\ge0}L_i$, with $L_0\simeq\FF$, $L_i L_j\subset L_{i+j}$ for all $i,j$, and $L_l=(L_1)^l$ for all $l\ge 1$. Choose
$$
\Pi:=\bigoplus_{i=1}^{\nu-1}L_i,\quad L:=L_1.
$$
In this case for any choice of the basis $e_1,\dots,e_k$ in $L$ the ideal $I$ is homogeneous, i.e.~generated by homogeneous relations. For an arbitrary nil-polynomial $P$ on $\Pi$ its restriction to $L$ expressed in the corresponding coordinates is
$$
R(y_1,\dots,y_k)=\frac{1}{\nu!}\omega\Bigl((y_1e_1+\dots+y_ke_k)^{\nu}\Bigr).
$$
Identifying $L$ with ${\mathfrak m}/{\mathfrak m}^2$, we observe that the homogeneous inverse system $R$ is a coordinate representation of the form ${\bf P}^{[\nu]}$ and therefore that of the form $a$ up to scale.

Theorem \ref{main} yields the well-known fact (see, e.g.~Proposition 7 in \cite{Em}) that a standard graded Artinian Gorenstein algebra, when written as a quotient by a homogeneous ideal, admits a homogeneous inverse system (note that any two such systems are proportional). In this situation the criterion for isomorphism of algebras in terms of inverse systems stated in Proposition 16 in \cite{Em} becomes rather simple: two quotients are isomorphic if and only if their homogeneous inverse systems are linearly equivalent up to scale (see Proposition 17 in \cite{Em} and Proposition 2.2 in \cite{ER}). We note that this classical criterion can be easily derived from Theorem \ref{main} as well. 
\end{remark} 

\section{The system of invariants}\label{sect2}
\setcounter{equation}{0}

In this section we survey the construction and properties of the system of invariants introduced in \cite{EI}. Here we assume that $\FF=\CC$, although much of what follows works for any algebraically closed field of characteristic zero.

Let $A$ and $\tilde A$ be graded Artinian Gorenstein algebras of vector space dimension greater than 2, and $P:\Pi\ra\CC$, $\tilde P:\tilde\Pi\ra\CC$ some nil-polynomials arising from $A$, $\tilde A$, respectively. Assume that $A$ and $\tilde A$ are isomorphic. As stated in Section \ref{sect1}, in this case $P$ and $\tilde P$ are linearly equivalent up to scale, i.e.~there exist $c\in\CC^*$ and a linear isomorphism $\varphi:\Pi\ra\tilde\Pi$ with $cP=\tilde P\circ\varphi$. Moreover, as shown in \cite{FIKK}, \cite{FK}, \cite{I}, the map
$$
\hat\varphi:{\mathfrak m}\ra\tilde{\mathfrak m},\quad x+y\mapsto\varphi(x)+c\,\tilde\omega_0^{-1}(\omega(y)),\quad x\in\Pi,\, y\in\Ann({\mathfrak m}),
$$
is an algebra isomorphism, where $\tilde\omega_0:=\tilde\omega|_{\Ann(\tilde{\mathfrak m})}$.

Further, for $s=2,\dots,\nu$ consider the forms ${\bf P}^{[s]}$ and $\tilde{\bf P}^{[s]}$ on ${\mathfrak m}/{\mathfrak m}^{\nu+2-s}$ and $\tilde{\mathfrak m}/\tilde{\mathfrak m}^{\nu+2-s}$ arising from $P$ and $\tilde P$, respectively, as explained in Section \ref{sect1}. Since the map $\hat\varphi$ is an algebra isomorphism, there exist algebra isomorphisms $\varphi^{[s]}: {\mathfrak m}/{\mathfrak m}^{\nu+2-s}\ra \tilde{\mathfrak m}/\tilde{\mathfrak m}^{\nu+2-s}$ such that $c{\bf P}^{[s]}={\bf P}^{[s]}\circ \varphi^{[s]}$. This fact allows one to utilize classical invariant theory for constructing numerical invariants of graded Artinian Gorenstein algebras. We will now recall the definitions of relative and absolute classical invariants (see, e.g.~\cite{O} for details). These definitions can be given in a coordinate-free setting. 

Let $W$ be a finite-dimensional complex vector space and ${\mathcal Q}_W^{m}$ the linear space of holomorphic forms of degree $m$ on $W$, with $m\ge 2$. Define an action of $\GL(W)$ on ${\mathcal Q}_W^{m}$ by the formula
$$
(C,Q)\mapsto Q_C,\,\,  Q_C(w):=Q(C^{-1}w),\,\,\hbox{where $C\in\GL(W)$, $Q\in{\mathcal Q}_W^{m}$, $w\in W$.}
$$
If two forms lie in the same $\GL(W)$-orbit, they are called linearly equivalent. An invariant (or relative classical invariant) of forms of degree $m$ on $W$ is a polynomial ${\mathcal I}:{\mathcal Q}_W^{m}\ra\CC$ such that for any $Q\in{\mathcal Q}_W^{m}$ and any $C\in\GL(W)$ one has ${\mathcal I}(Q)=(\det C)^{\ell}{\mathcal I}(Q_C)$, where $\ell$ is a non-negative integer called the weight of ${\mathcal I}$. It follows that ${\mathcal I}$ is in fact homogeneous of degree $\ell\dim_{\CC}W/m$. Finite sums of relative invariants comprise the algebra of polynomial $\SL(W)$-invariants of ${\mathcal Q}_W^{m}$, called the algebra of invariants (or algebra of classical invariants) of forms of degree $m$ on $W$. By the Hilbert Basis Theorem, this algebra is finitely generated. For any two invariants ${\mathcal I}$ and ${\mathcal J}$, with ${\mathcal J}\not\equiv 0$, the ratio ${\mathcal I}/{\mathcal J}$ yields a rational function on ${\mathcal Q}_W^{m}$ that is defined, in particular, at the points where ${\mathcal J}$ does not vanish. If ${\mathcal I}$ and ${\mathcal J}$ have equal weights, this function does not change under the action of $\GL(W)$, and we say that ${\mathcal I}/{\mathcal J}$ is an absolute invariant (or absolute classical invariant) of forms of degree $m$ on $W$. 

If one fixes coordinates $z_1,\dots,z_n$ in $W$, then $W$ is identified with $\CC^n$, $\GL(W)$ with $\GL(n,\CC)$, and any element $Q\in{\mathcal Q}^{m}_W$ is written as a homogeneous polynomial of degree $m$ in $z_1,\dots,z_n$. Invariants are usually defined in terms of the coefficients of the polynomial in $z_1,\dots,z_n$ representing $Q$. Observe, however, that the value of any absolute invariant at $Q$ is independent of the choice of coordinates in $W$.

The above discussion yields the following result. 

\begin{theorem}\label{theorem1}{\rm \cite{EI}} \sl Let $A$ be a graded Artinian Gorenstein algebra with $\dim_{\CC}A>2$, and $P$ a nil-polynomial arising from $A$. Further, for a fixed $s\in\{2,\dots,\nu\}$, let $W$ be a complex vector space isomorphic to ${\mathfrak m}/{\mathfrak m}^{\nu+2-s}$ by means of a linear map $\psi: W\ra {\mathfrak m}/{\mathfrak m}^{\nu+2-s}$. Fix an absolute invariant ${\mathbf I}$ of forms of degree $s$ on $W$. Then the value ${\mathbf I}(\psi^*{\bf P}^{[s]})$ depends only on the isomorphism class of $A$.  
\end{theorem}
For each $s$, let $\psi_s:\CC^{N_s}\ra{\mathfrak m}/{\mathfrak m}^{\nu+2-s}$ be some linear isomorphism, where $N_s:=\dim_{\CC}{\mathfrak m}/{\mathfrak m}^{\nu+2-s}$. We have thus constructed the following system of invariants:
$$
\displaystyle{\mathfrak I}:=\bigsqcup_{s=2}^{\nu}{\mathfrak I}_s,
$$
where
$$
\displaystyle{\mathfrak I}_s:=\Bigl\{{\mathbf I}(\psi_s^*{\bf P}^{[s]}),\,\hbox{${\bf I}$ is an absolute invariant of forms of degree $s$ on $\CC^{N_s}$}\Bigr\}.
$$ 

The highest stratum ${\mathfrak I}_{\nu}$ of ${\mathfrak I}$ was introduced in \cite{Ea} by means of the form $a$, which, as we noted in Section \ref{sect1}, coincides with any ${\bf P}^{[\nu]}$ up to scale. For standard graded Artinian Gorenstein algebras the idea of building invariants by the above method was briefly indicated in \cite{Em} in relation to homogeneous inverse systems. Any homogeneous inverse system in $\CC[z_1,\dots,z_k]$ arising from a given algebra of embedding dimension $k$ is calculated as explained in Remark \ref{stangraded} and therefore is proportional to a coordinate representation of any ${\bf P}^{[\nu]}$. Hence, the invariants resulting from the idea expressed in \cite{Em} also comprise the stratum ${\mathfrak I}_{\nu}$. This stratum will play an important role below. 

We now return to quasihomogeneous singularities, which are the main motivation of this study, and denote by ${\mathfrak I}^M$ the restriction of the system ${\mathfrak I}$ to the Milnor algebras of such singularities. In the remainder of the paper we will discuss the completeness property of ${\mathfrak I}^M$. At this stage, all our completeness results only concern homogeneous singularities, and this is the case that we will consider from now on. We thus further restrict ${\mathfrak I}^M$ to the class of Milnor algebras of homogeneous singularities and denote the restriction by ${\mathfrak I}^{MH}$.  

Let $Q$ be a non-zero element of ${\mathcal Q}_{\CC^n}^m$ with $n\ge 2$, $m\ge 3$, and $V_Q$ the germ at the origin of the hypersurface $\{Q=0\}$. Then the singularity of $V_Q$ is isolated if and only if $\Delta(Q)\ne 0$, where $\Delta$ is the discriminant (see Chapter 13 in \cite{GKZ}). Define
$$
X_n^m:=\{Q\in{\mathcal Q}_{\CC^n}^m:\Delta(Q)\ne 0\}.
$$
Any hypersurface germ $V$ at the origin in $\CC^n$ with homogeneous singularity and $\dim_{\CC}M(V)>1$ is biholomorphic to some $V_Q$ with $Q\in X_n^m$, $m\ge 3$. 

Next, for $Q,\tilde Q\in X_n^m$ the germs $V_Q$, $V_{\tilde Q}$ are biholomorphically equivalent if and only if $Q$ and $\tilde Q$ are linearly equivalent. On the other hand, we have the following fact.
\begin{proposition}\label{equivarbitr} {\rm \cite{EI}}\footnote{The proof of this proposition given in \cite{EI} was suggested to us by A. Gorinov.} 
\sl The orbits of the $\GL(n,\CC)$-action on $X_n^m$ are separated by invariant regular functions on the affine algebraic variety $X_n^m$, i.e.~by absolute invariants of the form $I/\Delta^{p}$, where $p$ is a non-negative integer and $I$ a relative invariant.
\end{proposition}
Let ${\mathcal I}_n^m$ be the algebra of invariant regular functions on $X_n^m$. By the Hilbert Basis Theorem, this algebra is finitely generated. The above discussion yields that the completeness of the system ${\mathfrak I}^{MH}$ will follow if one shows that the algebra ${\mathcal I}_n^m$ can be somehow \lq\lq extracted\rq\rq\, from ${\mathfrak I}^{MH}$.

Observe that $\emb\dim M(V_Q)=n$ and $M(V_Q)$ is canonically isomorphic to the quotient ${\mathbb M}(V_Q):=\CC[z_1,\dots,z_n]/{\mathbb J}(Q)$, where ${\mathbb J}(Q)$ is the ideal in $\CC[z_1,\dots,z_n]$ generated by all first-order partial derivatives of $Q$. From now on, we will only consider algebras of this form. Notice that the ideal ${\mathbb J}(Q)$ is homogeneous, and therefore, as stated in Remark \ref{stangraded}, the isomorphism class of ${\mathbb M}(V_Q)$ is determined by the linear equivalence class of any of its homogeneous inverse systems. Hence, it is a reasonable idea to explore whether ${\mathcal I}_n^m$ can be derived from the highest stratum ${\mathfrak I}_{\nu}^{MH}$ of ${\mathfrak I}^{MH}$.

By Lemma 3.3 in \cite{Sa}, the annihilator $\Ann({\mathfrak m})$ of the maximal ideal ${\mathfrak m}$ of ${\mathbb M}(V_Q)$ is generated by the element represented by the Hessian of $Q$, which implies that the nil-index of ${\mathfrak m}$ is found from the formula $\nu=n(m-2)$. Therefore, every nil-polynomial $P$ arising from ${\mathbb M}(V_Q)$ has degree $n(m-2)$, and, since  $\emb\dim{\mathbb M}(V_Q)=n$, the corresponding highest-degree form ${\bf P}^{[n(m-2)]}$ is a form on an $n$-dimensional vector space. We say that any such form ${\bf P}^{[n(m-2)]}$ is {\it associated}\, to $Q$ (recall that all these forms are proportional to each other). Thus, upon identification of ${\mathfrak m}/{\mathfrak m}^2$ with $\CC^n$, every element of ${\mathfrak I}_{n(m-2)}^{MH}$ calculated for the algebra ${\mathbb M}(V_Q)$ is given as ${\mathbf I}({\mathbf  Q})$, where ${\mathbf I}$ is an absolute invariant of forms of degree $n(m-2)$ on $\CC^n$ and ${\mathbf Q}$ is any form associated to $Q$.

Computing associated forms is quite easy for any realization of the algebra. Choose a basis $e_1,\dots,e_n$ in a complement to ${\mathfrak m}^2$ in ${\mathfrak m}$, and let $f_j$ be the element of ${\mathfrak m}/{\mathfrak m}^2$ represented by $e_j$ for $j=1,\dots,n$. Denote by $w_1,\dots,w_n$ the coordinates in ${\mathfrak m}/{\mathfrak m}^2$ with respect to the basis $f_1,\dots,f_n$. Further, choose a vector $v$ spanning $\Ann({\mathfrak m})$. If $k_1,\dots,k_n$ are non-negative integers such  that $k_1+\dots+k_n=n(m-2)$, the product $e_1^{k_1}\dots e_n^{k_n}$ is an element of $\Ann({\mathfrak m})$, and thus we have $e_1^{k_1}\dots e_n^{k_n}=\mu_{k_1,\dots,k_n}v$ for some $\mu_{k_1,\dots,k_n}\in\CC$. Then the form
\begin{equation}
\sum_{k_1+\dots+k_n=n(m-2)}\mu_{k_1,\dots,k_n}\left(
\begin{array}{c}
n(m-2)\\
k_1,\dots,k_n
\end{array}
\right)
w_1^{k_1}\dots w_n^{k_n}\label{assoccomp}
\end{equation}
is a coordinate representation of a form associated to $Q$, where
$$
\left(
\begin{array}{c}
n(m-2)\\
k_1,\dots,k_n
\end{array}
\right):=\displaystyle\frac{(n(m-2))!}{k_1!\dots k_n!}
$$
is a multinomial coefficient.

We now propose a conjecture.

\begin{conjecture}\label{conj1} \rm For any ${\tt I}\in{\mathcal I}_n^m$ there exists an absolute invariant ${\mathbf I}$ of forms of degree $n(m-2)$ on $\CC^n$ such that for all $Q\in X_n^m$ the invariant ${\mathbf I}$ is defined at some (hence at every) form ${\bf Q}$ associated to $Q$ and ${\mathbf I}({\mathbf Q})={\tt I}(Q)$.
\end{conjecture}
For binary quartics ($n=2$, $m=4$) and ternary cubics ($n=3$, $m=3$) the conjecture was essentially verified in \cite{Ea} (see also \cite{EI} and Example \ref{examses} below). Furthermore, in \cite{EI} we showed that the conjecture holds for binary quintics ($n=2$, $m=5$) and binary sextics ($n=2$, $m=6$) as well. If the conjecture were confirmed in full generality, it would imply that the stratum ${\mathfrak I}_{n(m-2)}^{MH}$ is a complete system of invariants for homogeneous hypersurface singularities defined by forms in $X_n^m$.

We will now mention another interesting consequence of Conjecture \ref{conj1}. First, we observe that ${\mathbf I}({\mathbf Q})$ is rational when regarded as a function of $Q$. In order to see this, we choose a particular (canonical) identification of ${\mathfrak m}/{\mathfrak m}^2$ with $\CC^n$. Namely, for $j=1,\dots,n$ let $e_j$ be the element of ${\mathfrak m}$ represented by the coordinate function $z_j$. Also, we let $v$ be the element represented by the Hessian of $Q$. Denote by ${\bf Q}^c$ the corresponding associated form found from formula (\ref{assoccomp}). As follows from Remark \ref{stangraded}, the form ${\bf Q}^c$ is an inverse system for ${\mathbb M}(V_Q)$. It is clear that $\mu_{k_1,\dots,k_n}$ that occur in (\ref{assoccomp}) are rational functions of the coefficients of $Q$, and therefore for an absolute invariant ${\mathbf I}$ of forms of degree $n(m-2)$ on $\CC^n$ the expression ${\mathbf I}({\mathbf Q}^c)$ is a rational function of $Q$. 

Let ${\mathcal R}_n^m$ denote the collection of all invariant rational functions on $X_n^m$ obtained in this way. Further, let $\hat {\mathcal I}_n^m$ be the algebra of restrictions to $X_n^m$ of all absolute invariants of forms of degree $m$ on $\CC^n$. Note that ${\mathcal R}_n^m$ lies in $\hat{\mathcal I}_n^m$ (see Proposition 1 in \cite{DC}). We claim that Conjecture \ref{conj1} implies
\begin{equation}
{\mathcal R}_n^m=\hat {\mathcal I}_n^m.\label{conjweaker}
\end{equation} 
Indeed, since every element of $\hat {\mathcal I}_n^m$ can be represented as a ratio of two elements of ${\mathcal I}_n^m$ (see Proposition 6.2 in \cite{Mu}), identity (\ref{conjweaker}) is equivalent to the inclusion ${\mathcal I}_n^m\subset{\mathcal R}_n^m$, which clearly follows from Conjecture \ref{conj1} (cf. Conjecture 3.2 in \cite{EI}). 

We remark that identity (\ref{conjweaker}) is interesting from the invariant-theoretic point of view, since it means that the invariant theory of forms of degree $m$ can be completely recovered from that of forms of degree $n(m-2)$. Observe that (\ref{conjweaker}) is {\it a priori}\, weaker than Conjecture \ref{conj1}. Indeed, it may potentially happen that for some ${\tt I}\in{\mathcal I}_n^m$ there exist an absolute invariant ${\mathbf I}$ of forms of degree $n(m-2)$ on $\CC^n$ such that ${\mathbf I}({\mathbf Q})\equiv{\tt I}(Q)$, where ${\bf I}({\mathbf Q})$ is regarded as a function of $Q$, but for some $Q_0\in X_n^m$ the invariant ${\bf I}$ is not defined at the forms associated to $Q_0$. In a forthcoming paper by J. Alper and the second author identity (\ref{conjweaker}) will be shown to hold for binary forms of any degree.

We will now illustrate Conjecture \ref{conj1} with the example of simple elliptic singularities of type $\tilde E_6$. 

\begin{example}\label{examses}\rm Simple elliptic $\tilde E_6$-singularities form a family $V_t$ parametrized by $t\in\CC$ satisfying $t^3+27\ne 0$. Namely, for every such $t$ let $V_t:=V_{Q_{{}_t}}$, where $Q_t$ is the following cubic on $\CC^3$:
$$
Q_t(z_1,z_2,z_3):=z_1^3+z_2^3+z_3^3+tz_1z_2z_3.
$$
Since $n=m=3$, we have $\nu=n(m-2)=3$, thus any form associated to $Q_t$ is again a ternary cubic. To compute such a form using formula (\ref{assoccomp}), set $e_j$ to be the element of ${\mathfrak m}$ represented by $z_j$ for $j=1,2,3$ and $v$ the element represented by $z_1z_2z_3$. Then for the coefficients in formula (\ref{assoccomp}) we have
$$
\mu_{3,0,0}=\mu_{0,3,0}=\mu_{0,0,3}=-\frac{t}{3},\quad\mu_{1,1,1}=1,
$$
with all the remaining $\mu_{k_1,k_2,k_3}$ being zero. These coefficients yield the following associated form:
$$
{\bf Q}_t:=-\frac{t}{3}(w_1^3+w_2^3+w_3^3)+6w_1w_2w_3.
$$
The form ${\bf Q}_t$ is an inverse system for ${\mathbb M}(V_t)$ and has been known for a long time (see \cite{Ea}, \cite{Em}). For $t\ne 0$, $t^3-216\ne 0$ one has $\Delta({\bf Q}_t)\ne 0$, in which case the original cubic $Q_t$ is associated to ${\bf Q}_t$ and thus there is a natural duality between $Q_t$ and ${\bf Q}_t$.

Further, the algebra of classical invariants of ternary cubics is generated by certain invariants ${\mathcal I}_4$ and ${\mathcal I}_6$, where the subscripts indicate the degrees (see pp.~381--389 in \cite{Ell}). For a ternary cubic of the form
$$
Q(z)=az_1^3+bz_2^3+cz_3^3+6dz_1z_2z_3
$$
the values of ${\mathcal I}_4$ and ${\mathcal I}_6$ are computed as follows: 
$$
\begin{array}{l}
{\mathcal I}_4(Q)=abcd-d^4 ,\quad {\mathcal I}_6(Q)=a^2b^2c^2-20abcd^3-8d^6 ,
\end{array}
$$
and $\Delta(Q)={\mathcal I}_6^2+64{\mathcal I}_4^3$.\footnote{This formula for the discriminant of a ternary cubic differs from the general formula given in \cite{GKZ} by a scalar factor.} 

Consider the $j$-invariant of ternary cubics, which is the absolute invariant defined as follows:
$$
j:=\frac{64\,{\mathcal I}_4^3}{\Delta}.\label{invariantj}
$$
It is easy to see that the restriction $j|_{X_3^3}$ generates the algebra ${\mathcal I}_3^3$. In particular, any two non-equivalent ternary cubics with non-vanishing discriminant are distinguished by $j$ (see Proposition \ref{equivarbitr}). Further, we have
$$
j(Q_t)=-\frac{t^3(t^3-216)^3}{1728(t^3+27)^3}.
$$
Details on computing $j(Q)$ for any ternary cubic $Q$ with $\Delta(Q)\ne 0$ can be found, for example, in \cite{Ea}.  

Next, consider the following absolute invariant of ternary cubics:
$$
{\bf j}:=\frac{1}{j}.
$$     
A straightforward calculation shows that for any $Q\in X_3^3$ the absolute invariant ${\bf j}$ is defined at ${\bf Q}_t$ and ${\bf j}({\bf Q}_t)=j(Q_t)$, which demonstrates that Conjecture \ref{conj1} is indeed valid for $n=m=3$. 
\end{example}

{\obeylines
\noindent Mathematical Sciences Institute
\noindent The Australian National University
\noindent Canberra, ACT 0200
\noindent Australia
\noindent e-mail: michael.eastwood@anu.edu.au, alexander.isaev@anu.edu.au
}


\begin{thebibliography}{ABCD}

\bibitem[B]{B} Bass, H., On the ubiquity of Gorenstein rings, {\it Math. Z.} 82 (1963), 8--28.

\bibitem[BH]{BH} Bruns, W. and Herzog, J., {\it Cohen-Macaulay Rings}, Cambridge Studies in Advanced Mathematics 39, Cambridge University Press, Cambridge, 1993.

\bibitem[DC]{DC} Dieudonn\'e, J. A. and Carrell, J. B., Invariant theory, old and new, {\it Adv. in Math.}  4 (1970), 1--80.

\bibitem[Ea]{Ea} Eastwood, M. G., Moduli of isolated hypersurface singularities, {\it Asian J. Math.} 8 (2004), 305--313.

\bibitem[EI]{EI} Eastwood, M. G. and Isaev, A. V., Extracting invariants of isolated hypersurface singularities from their moduli algebras, to appear in {\it Math. Ann.}, DOI 10.1007/s00208-012-0836-7.

\bibitem[ER]{ER} Elias, J. and Rossi, M. E., Isomorphism classes of short Gorenstein local rings via Macaulay's inverse system, {\it Trans. Amer. Math. Soc.} 364 (2012), 4589--4604.

\bibitem[Ell]{Ell} Elliott, E. B., {\it An Introduction to the Algebra of Quantics}, Oxford University Press, 1895.

\bibitem[Em]{Em} Emsalem, J., G\'eom\'etrie des points \'epais, {\it Bull. Soc. Math. France} 106 (1978), 399--416.

\bibitem[FIKK]{FIKK} Fels, G., Isaev, A., Kaup, W. and Kruzhilin, N., Isolated hypersurface singularities and special polynomial realizations of affine quadrics, {\it J. Geom. Analysis} 21 (2011), 767--782.

\bibitem[FK]{FK} Fels, G. and Kaup, W., Nilpotent algebras and affinely homogeneous surfaces, {\it Math. Ann.} 353 (2012), 1315--1350.

\bibitem[GKZ]{GKZ} Gelfand, I. M., Kapranov, M. M. and Zelevinsky, A. V., {\it Discriminants, Resultants and Multidimensional Determinants}, Modern Birkh\"auser Classics, Birkh\"auser Boston, Inc., Boston, MA, 2008.

\bibitem[GLS]{GLS} Greuel, G.-M., Lossen, C. and Shustin, E., {\it Introduction to Singularities and Deformations}, Springer Monographs in Mathematics, Springer, Berlin, 2007.

\bibitem[He]{He} Hertling, C., {\it Frobenius Manifolds and Moduli Spaces for Singularities}, Cambridge Tracts in Mathematics 151, Cambridge University Press, Cambridge, 2002.

\bibitem[Hu]{Hu} Huneke, C., Hyman Bass and ubiquity: Gorenstein rings, in {\it Algebra, $K$-theory, Groups, and Education} (New York, 1997), Contemp. Math. 243, Amer. Math. Soc., Providence, RI, 1999, pp. 55--78.

\bibitem[I]{I} Isaev, A. V., On the affine homogeneity of algebraic hypersurfaces arising from Gorenstein algebras, {\it Asian J. Math.} 15 (2011), 631--640.


\bibitem[MY]{MY} Mather, J. and Yau, S. S.-T., Classification of isolated hypersurface singularities by their moduli algebras, {\it Invent. Math.} 69 (1982), 243--251.

\bibitem[Mu]{Mu} Mukai, S., {\it An Introduction to Invariants and Moduli}, Cambridge Studies in Advanced Mathematics 81, Cambridge University Press, Cambridge, 2003.

\bibitem[O]{O} Olver, P., {\it Classical Invariant Theory}, London Mathematical Society Student Texts 44, Cambridge University Press, Cambridge, 1999.

\bibitem[Sa]{Sa} Saito, K., Einfach-elliptische Singularit\"aten, {\it Invent. Math.} 23 (1974), 289--325.

\bibitem[SV]{SV} Sharpe, D. W. and V\'amos, P., {\it Injective Modules}, Cambridge Tracts in Mathematics and Mathematical Physics 62, Cambridge University Press, London--New York, 1972.

\bibitem[Sh]{Sh} Shoshitaishvili, A. N., Functions with isomorphic Jacobian ideals, {\it Funct. Anal. Appl.} 10 (1976), 128--133.

\end{thebibliography}
\end{document}